# The spectrum of the variety of anti-rectangular Abel Grassmann bands

R. A. R. Monzo
Flat 10 Albert Mansions
Crouch Hill, London N8 9RE
United Kingdom

Abstract

We prove that the set of all orders of finite algebras in the groupoid variety of anti-rectangular Abel Grassmann bands consists of all powers of four. We also prove that any groupoid anti-isomorphic to a finite or countable anti-rectangular Abel Grassmann band $G$ is isomorphic to $G$. It is proved that within isomorphism there is only one countable anti-rectangular Abel Grassmann band and that it is isomorphic to a proper subset of itself.

**1. Introduction.**

Kazim and Naseeruddin studied a groupoid variety consisting of what they called *"almost semigroups"*, groupoids satisfying the equation $(xy)z = (zy)x$ [8]. Such groupoids have also been referred to as *"invertive"*, *"Abel-Grassmann's"* and *"right modular"* and *"left almost semigroups"* [4,16,17,6,9]. Various aspects of these $AG$-groupoids have been studied over the years, such as partial ordering and congruences, inflations, bands, zeroids and idempoids, ideals, topological structure, various forms of regularity, maximal separative homomorphic images and the least semilattice congruence, fuzziness, structure of unions of groups, power groupoids and inclusion classes and simplicity [1,5,6,7,9-17]. In this paper we study $AG$-groupoids that satisfy the additional equations $x = x^2$ and $(xy)x = y$. These groupoids are called *"anti-rectangular $AG$-bands"* [15].

The set of orders of the finite algebras in a groupoid variety $V$ is called the **spectrum of $V$**. We will denote this by $sp(V)$. T. Evans showed that the spectrum of the groupoid variety defined by the equation $(xy)(yz) = y$ is the set $\{1, 4, 9, 16, ...\} = \{n^2\}_{n \in N}$ of all squares [3]. Evans generalised this result and obtained, for each positive integer $n \in N = \{1, 2, 3, ...\}$, a variety of groupoids having as spectrum all $n^{th}$ powers [2]. Herein we determine the spectrum of the groupoid variety determined by the equations $(xy)z = (zy)x$, $x = x^2$ and $(xy)x = y$, the variety of anti-rectangular Abel Grassmann bands, which we will denote by $ARAGB$. We prove below that $sp(ARAGB) = \{4^n\}_{n \in N \cup \{0\}}$. Comparing this to Evans' results one can wonder whether, for any $m^{th}$ power $q^m$, there is a groupoid variety $V$ with $sp(V) = \{(q^m)^n\}_{n \in N \cup \{0\}}$. This question will <u>not</u> be dealt with in this paper.

There is another reason to study the structure of anti-rectangular Abel Grassmann bands. Let us denote the variety of Abel Grassmann groupoids as $AG$. That is, $AG$ is the variety of groupoids determined by the equation $(xy)z = (zy)x$. Protic and Stepanovic proved that any $AG$-band $B$ is an $AG$-band $Y$ of anti-rectangular $AG$-bands [15]. That is,

**Result 1. [15, Th. 2.1]** *If $B \in AG$ and $x = x^2$ $(x \in B)$ then there exists a band $Y \in AG$ such that $B = \overset{\bullet}{\cup} G_\alpha (\alpha \in Y)$, $G_\alpha G_\beta \subseteq G_{\alpha\beta} (\alpha, \beta \in Y)$, and $G_\alpha \in ARAGB$ $(\alpha \in Y)$.*

Our results will imply that there is, up to isomorphism, exactly one anti-rectangular $AG$-band of order $4^n$ for each $n \in \{0, 1, 2, ...\}$ and that there are no finite, anti-rectangular $AG$-bands of any other orders.

So the finite, anti-rectangular $AG$-bands are a basic building block of the finite $AG$-bands. As we shall see the basic building block of the finite, anti-rectangular $AG$-bands is the following anti-rectangular $AG$-band $G$ of order 4. It is isomorphic to the anti-rectangular $AG$-band generated by any two distinct elements, $a$ and $b$ say, of any anti-rectangular $AG$-band:

**Example 1.**

| G | a | b | ab | ba |
|---|---|---|----|----|
| a | a | ab | ba | b |
| b | ba | b | a | ab |
| ab | b | ba | ab | a |
| ba | ab | a | b | ba |

We will also show that if $H \in \mathbf{ARAGB}$ and $|H| = 4^n$ then $H$ consists of exactly $4^{n-1}$ disjoint copies of $G$. The following results will be used throughout the paper.

**Result 2.** *If $R \in \mathbf{AG}$ and $\{c,d,e,f\} \subseteq R$ then $(cd)(ef) = (ce)(df)$.*

**Result 3.** *If $R \in \mathbf{ARAGB}$ and $\{a,b,c\} \subseteq R$ then $a(bc) = c(ba)$.*

*Proof*: $c(ba) = [(ba)c]c = [(ca)b]c = \{(ca)[(cb)c]\}c =$
$= \{[c(cb)](ac)\}c = [c(ac)][c(cb)] = a(bc)$ ∎

**Result 4.** *Let $R \in \mathbf{ARAGB}$, with $\{c,d\} \subseteq R$ and $c \neq d$. Then the subgroupoid $\langle c,d \rangle$ of R generated by c and d is isomorphic to the groupoid G in Example 1. One isomorphism is given by the mapping $c \to a$, $d \to b$, $cd \to ab$ and $dc \to ba$.*

**Result 5.** *Any two distinct elements of G generate G.*

**Result 6.** *Any bijection on G is either an isomorphism or an anti-isomorphism. Four-cycles and two-cycles are anti-isomorphisms and the identity mapping, three-cycles and products of two-cycles are isomorphisms.*

**Result 7.** *Any groupoid anti-isomorphic to G is isomorphic to G. In particular, if $\Phi: G \to H$ is an anti-isomorphism then the mapping $a \to \Phi a$, $b \to \Phi b$, $ab \to \Phi(ba)$ and $ba \to \Phi(ab)$ is an isomorphism.*

**Result 8.** *Suppose that A and B are subgroupoids of $R \in \mathbf{ARAGB}$ and that A and B are isomorphic to G. Then either $A = B$, $A \cap B = \emptyset$ or $A \cap B = \{c\}$.*

**Notation:** $R \cong R'$ [ $R \cong \overline{R'}$ ] will denote that $R$ and $R'$ are [anti-] isomorphic.

**Result 9.** *If $R \in \mathbf{ARAGB}$ and $R \cong \overline{R}$ then $\overline{R} \in \mathbf{ARAGB}$.*

*Proof*: Let $\Phi: R \to \overline{R}$ be an anti-isomorphism. Then it is straightforward to show that $\overline{R}$ is a band that satisfies the equation $(xy)x = y$. Let $\{r_1', r_2', r_3'\} \subseteq \overline{R}$. Then there exists $\{r_1, r_2, r_3\} \subseteq R$ such that $r_i' = \Phi r_i \, (i \in \{1,2,3\})$. Using Result 3,
$(r_1' r_2') r_3' = [(\Phi r_1)(\Phi r_2)](\Phi r_3) = [\Phi(r_2 r_1)](\Phi r_3) = \Phi[r_3(r_2 r_1)] = \Phi[r_1(r_2 r_3)] =$
$= [\Phi(r_2 r_3)](\Phi r_1) = [(\Phi r_3)(\Phi r_2)](\Phi r_1) = (r_3' r_2') r_1'$ and so $\overline{R}$ satisfies the equation $(xy)z = (zy)x$. Hence, $\overline{R} \in ARAGB$. ∎

**2. Finite anti-rectangular Abel Grassmann bands.**

We use $H \leq R \, (H \prec R)$ to denote that H is a subgroupoid (proper subgroupoid) of the groupoid *R*. Suppose that $G \leq H \prec R$ and $\{H, R\} \subseteq ARAGB$. Let $x \in R - H$. If $H_x = H \cup \{xh\}_{h \in H} \cup \{hx\}_{h \in H} \cup \{(ax)h\}_{h \in H}$ then $H_x \in ARAGB$, as follows from the following multiplication table:

**Table 1.** (for $\{h, \overline{h}\} \subseteq H$)

| $H_x$ | $\overline{h}$ | $x\overline{h}$ | $\overline{h}x$ | $(ax)\overline{h}$ |
|---|---|---|---|---|
| $h$ | $h\overline{h}$ | $(ax)[(\overline{h}a)h]$ | $x(\overline{h}h)$ | $[(h\overline{h})(ah)]x$ |
| $xh$ | $(\overline{h}h)x$ | $x(h\overline{h})$ | $(ax)[\overline{h}(ah)]$ | $a(h\overline{h})$ |
| $hx$ | $(ax)[(ha)(\overline{h}h)]$ | $\overline{h}h$ | $(h\overline{h})x$ | $x[(ah)\overline{h}]$ |
| $(ax)h$ | $x[h(\overline{h}a)]$ | $[(h\overline{h})a]x$ | $(ah)(\overline{h}a)$ | $(ax)(h\overline{h})$ |

*Proof*: We will only calculate the products in the first row. Calculations for products in the other three rows are along similar lines, using Result 1 together with the fact that $R \in ARAGB$. Firstly,

$h(x\overline{h}) = [(x\overline{h})h]h = [(h\overline{h})x]h = [h\overline{h}h](xh) = \overline{h}(xh) = [a(\overline{h}a)](xh) = (ax)[(\overline{h}a)h]$

Also, $h(\overline{h}x) = [(\overline{h}x)h]h = [(hx)\overline{h}]h = (hxh)(\overline{h}h) = x(\overline{h}h)$. Finally,

$(h\overline{h})(ah) = (ha)(\overline{h}h) = [h(\overline{h}h)][a(\overline{h}h)] = \overline{h}[a(\overline{h}h)] = (\overline{h}a)(\overline{h}h) = [(h\overline{h})a]\overline{h} =$
$= [(h\overline{h})(x\overline{h})](ax) = [(hx)\overline{h}](ax) = [(ax)\overline{h}](hx)$ ; therefore,

$$\left[(h\bar{h})(ah)\right]x = \left\{\left[(ax)\bar{h}\right](hx)\right\}x = (xhx)\left[(ax)\bar{h}\right] = h\left[(ax)\bar{h}\right]. \blacksquare$$

Note that it is straightforward to show that the sets $H, \{xh\}_{h\in H}$, $\{hx\}_{h\in H}$ and $\{(ax)h\}_{h\in H}$ are pairwise disjoint sets. Furthermore, it is easy to show that for $\{h, h'\} \subseteq H$, two elements $xh$ and $xh'$ [$hx$ and $h'x$; $(ax)h$ and $(ax)h'$] are equal if and only if $h = h'$. Therefore, if $H$ contains $n$ elements then $H_x$ contains $4n$ elements. We will call $H_x$ **the extension of $H$ by $x$.** We have shown that:

**Theorem 1.** *If $G \leq H \prec R$ and $\{H, R\} \subseteq \mathbf{ARAGB}$, with $x \in R - H$ then $H_x$, the extension of $H$ by $x$, is in $\mathbf{ARAGB}$. Furthermore if $|H| = n$ then $|H_x| = 4n$.*

**Corollary 2.** $sp(\mathbf{ARAGB}) = \{4^n\}_{n \in N \cup \{0\}}$.

**Corollary 3.** *An anti-rectangular $\mathbf{AG}$-band of order $4^n$ has $(n+1)$-generators $(n \in \{0, 1, 2...\})$.*

**Theorem 4.** *Suppose that $G \leq H \in \mathbf{ARAGB}$ and $x \notin H$. We define pairwise disjoint sets $A = \{xh\}_{h\in H}$, $B = \{hx\}_{h\in H}$ and $C = \{(ax)h\}_{h\in H}$ such that $A \cap H = B \cap H = C \cap H = \emptyset$. Define $H^x = H \cup A \cup B \cup C$ with a product defined as in Table 1. Then $H^x \in \mathbf{ARAGB}$ and $H^x \cong H_x$.*

*Proof*: The product is well defined and so $H^x$ is a groupoid. To show that $H^x \in \mathbf{AG}$ we need to show that the products in Table 1 satisfy the equation $(yz)w = (wz)y$. This follows from the fact that $H \in \mathbf{AG}$. We show this for 2 of the 64 forms of products only:

The product $\left\{(x\bar{h})\left[(ax)\bar{\bar{h}}\right]\right\}(hx) = \left[a\left(\bar{h}\bar{\bar{h}}\right)\right](hx) = x\left\{h\left[a\left(\bar{h}\bar{\bar{h}}\right)\right]\right\}$; then $\left\{(hx)\left[(ax)\bar{\bar{h}}\right]\right\}(x\bar{h}) = \left\{x\left[(ah)\bar{\bar{h}}\right]\right\}(x\bar{h}) = x\left\{\left[(ah)\bar{\bar{h}}\right]\bar{h}\right\}$. Then, 

$$\left[(ah)\bar{\bar{h}}\right]\bar{h} = \left(\bar{h}\bar{\bar{h}}\right)(ah) = \left[\left(\bar{h}\bar{\bar{h}}\right)a\right]\left[\left(\bar{h}\bar{\bar{h}}\right)h\right] = \left\{\left[a\left(\bar{h}\bar{\bar{h}}\right)\right]\left(\bar{h}\bar{\bar{h}}\right)\right\}\left[\left(\bar{h}\bar{\bar{h}}\right)h\right] =$$
$$= \left\{\left[\left(\bar{h}\bar{\bar{h}}\right)h\right]\left(\bar{h}\bar{\bar{h}}\right)\right\}\left[a\left(\bar{h}\bar{\bar{h}}\right)\right] = h\left[a\left(\bar{h}\bar{\bar{h}}\right)\right] \text{ and so}$$

$$\{(x\bar{h})[(ax)\bar{\bar{h}}]\}(hx) = \{(hx)[(ax)\bar{\bar{h}}]\}(x\bar{h}).$$

Consider that $\{[(ax)h]\bar{h}\}(\bar{\bar{h}}x) = \{x[h(\bar{h}a)]\}(\bar{\bar{h}}x) = (ax)\{\bar{\bar{h}}[a\{h(\bar{h}a)\}]\}$ and that $[(\bar{\bar{h}}x)\bar{h}][(ax)h] = \{(ax)[(\bar{\bar{h}}a)(\bar{h}\bar{\bar{h}})]\}[(ax)h] = (ax)\{[(\bar{\bar{h}}a)(\bar{h}\bar{\bar{h}})]h\}$. Then
$\bar{\bar{h}}\{a[h(\bar{h}a)]\} = \bar{\bar{h}}[(ah)\bar{h}] = \bar{\bar{h}}[(\bar{h}h)a] = [\bar{\bar{h}}(\bar{h}h)](\bar{h}a)$ and
$[(\bar{\bar{h}}a)(\bar{h}\bar{\bar{h}})]h = [h(\bar{h}\bar{\bar{h}})](\bar{\bar{h}}a)$. But then $h(\bar{h}\bar{\bar{h}}) = [(\bar{h}\bar{\bar{h}})h]h = [(\bar{h}h)(\bar{h}\bar{\bar{h}})]h =$
$= [h(\bar{h}\bar{\bar{h}})](\bar{h}h) = \bar{\bar{h}}(\bar{h}h)$ and so $\{[(ax)h]\bar{h}\}(\bar{\bar{h}}x) = [(\bar{\bar{h}}x)\bar{h}][(ax)h]$.

Using Table 1 and the fact that $H \in \boldsymbol{ARAGB}$ it is straightforward to prove that $H^x$ satisfies the equations $y = y^2$ and $y(zy) = z$ and so $H^x \in \boldsymbol{ARAGB}$. Finally, since the multiplication tables of $H^x$ and $H_x$ are exactly the same they are isomorphic groupoids. ∎

**Definition 1.** We define $G_0$ as the trivial groupoid, $G_1 = G$ and by induction, $G_n \equiv G_{n-1}^{x_{n-1}}$ $(n \geq 2)$, where $x_{n-1} \notin G_{n-1}$.

**Corollary 5.** *Any finite anti-rectangular* $\boldsymbol{AG}$ *- band is isomorphic to* $G_n$ *for some* $n \in \boldsymbol{N} \cup \{0\}$. *If* $H \in \boldsymbol{ARAGB}$ *and* $|H| = 4n$ *then* $H \cong G_n$.

**Corollary 6.** *For* $n \geq 1$, $G_n$ *is a Y band of anti-rectangular* $\boldsymbol{AG}$ *- bands* $G_\alpha (\alpha \in Y)$ *where* $G_\alpha \cong G_{n-1} (\alpha \in Y)$ *and* $Y \cong G$.

## 3. Countable anti-rectangular Abel Grassmann bands

In this section we show that, to within isomorphism, there is precisely one countable anti-rectangular Abel Grassmann band. This result will follow from the following construction of such a groupoid.

**Construction 1.** *Let* $H = \bigcup_{n=1}^{\infty} G_n$, *with the* $G_n$'s *as in Definition 1. Define a product* $*$ *on* $H$ *as follows. If* $\{u,v\} \subseteq H$ *with* $u \in G_{n_u} - G_{n_u-1}$ *and* $v \in G_{n_v} - G_{n_v-1}$ *then* $u * v$ *is defined as the product of $u$ and $v$ in* $G_{\max\{n_u, n_v\}}$.

Clearly $*$ is well defined. Since, by Theorem 4, every $G_n \in ARAGB$, and since $\max\{\max\{n_u, n_v\}, n_w\} = \max\{\max\{n_w, n_v\}, n_u\}$, it is straightforward to prove that $H \in ARAGB$. It is also clear that $|H| = \chi_0$. Hence,

**Theorem 7.** The groupoid $H$ in Construction 1 is a countable anti-rectangular $AG$ - band.

**Theorem 8.** *Any countable anti-rectangular $AG$ -band $K$ is isomorphic to $H$ in Construction 1.*

*Proof.* Let $K = \bigcup_{n=1}^{\infty}\{y_n\}$, with $y_i = y_j$ if and only if $i = j$. Define $K_0 = \varnothing$, $K_1 = \{y_1, y_2, y_1 y_2, y_2 y_1\}$ and $R_1 = K - K_1$. Define $K_2 = K_1^{y_{t_1}}$, where $t_1$ is the minimum of the subscripts of the $y_n$'s in $R_1$. Define $R_2 = K - K_2$ and $K_3 = K_2^{y_{t_2}}$, where $t_2$ is the minimum subscript of the $y_n$'s in $R_2$. In general, by induction we define $R_n = K - K_n$ and $K_{n+1} = K_n^{y_{t_n}}$, where $t_n$ is the minimum subscript of the $y_n$'s in $R_n$. Then every $y_n$ must eventually appear in some $K_t$ and therefore $K = \bigcup_{n=0}^{\infty} K_n$. By Result 2, $\{y_1, y_2, y_1 y_2, y_2 y_1\}$ is isomorphic to $G_1 = G$. Call this isomorphism $\Phi_1$. Note that $\Phi_1 y_1 = a$, $\Phi_1 y_2 = b$, $\Phi_1(y_1 y_2) = ab$ and $\Phi_1(y_2 y_1) = ba$. Now by induction, we define $\Phi_n : K_n \to G_n$ as follows. Firstly, $\Phi_n = \Phi_{n-1}$ on $K_{n-1}$. Then for $k \in K_n - K_{n-1}$ we define $\Phi_n(y_{t_{n-1}} k) = x_{t_{n-1}}(\Phi_{n-1} k)$, $\Phi_n(k y_{t_{n-1}}) = (\Phi_{n-1} k) x_{t_{n-1}}$ and $\Phi_n[(y_1 y_{t_{n-1}}) k] = [(\Phi_{n-1} y_1) x_{t_{n-1}}](\Phi_{n-1} k)$. We now prove by induction on $n$ that $\Phi_n : K_n \to H_n$ is an isomorphism. If $n = 1$ then clearly $\Phi_1$ is an isomorphism. Assume that for $1 \leq t \prec n$, $\Phi_t$ is an isomorphism. Then the fact that $\Phi_n$ is 1-1 and onto $G_n$ follows from the definition of $\Phi_n$ and the fact that $\Phi_{n-1}$ is 1-1 and onto $H_{n-1}$. The fact that $\Phi_n(\alpha\beta) = (\Phi_n \alpha)(\Phi_n \beta)$ for any $\{\alpha, \beta\} \subseteq K_n$ follows from the definition of product in $K_n$ and $G_n$ and the fact that $\Phi_{n-1}$ is an isomorphism. We shall prove this for two of the 16 possible types of products in $K_n$. The proofs for the other products are similar. Recall that if $\{k, \bar{k}\} \subseteq K_{n-1}$, products in $K_n$ are:

**Table 2.** (for $\{k, \bar{k}\} \subseteq K_{n-1}$)

| $K_n = K_{n-1}^{y_{t_{n-1}}}$ | $\bar{k}$ | $y_{t_{n-1}} \bar{k}$ | $\bar{k} y_{t_{n-1}}$ | $(y_1 y_{t_{n-1}}) \bar{k}$ |
|---|---|---|---|---|
| $k$ | $k\bar{k}$ | $(y_1 y_{t_{n-1}})[(\bar{k} y_1) k]$ | $y_{t_{n-1}}(\bar{k} k)$ | $[(k\bar{k})(y_1 k)] y_{t_{n-1}}$ |
| $y_{t_{n-1}} k$ | $(\bar{k} k) y_{t_{n-1}}$ | $y_{t_{n-1}}(k\bar{k})$ | $(y_1 y_1 y_{t_{n-1}})[\bar{k}(y_1 k)]$ | $y_1(k\bar{k})$ |
| $k y_{t_{n-1}}$ | $(y_1 y_{t_{n-1}})[(k y_1)(\bar{k} k)]$ | $\bar{k} k$ | $(k\bar{k}) y_{t_{n-1}}$ | $y_{t_{n-1}}[(y_1 k)\bar{k}]$ |
| $(y_1 y_{t_{n-1}}) k$ | $y_{t_{n-1}}[k(\bar{k} y_1)]$ | $[(k\bar{k}) y_1] y_{t_{n-1}}$ | $(y_1 k)(\bar{k} y_1)$ | $(y_1 y_{t_{n-1}})(k\bar{k})$ |

So $\Phi_n\left[(ky_{t_{n-1}})\bar{k}\right] = \Phi_n\{(y_1 y_{t_{n-1}})[(ky_1)(\bar{k} k)]\} = [(\Phi_{n-1} y_1) x_{t_{n-1}}] \Phi_{n-1}[(ky_1)(\bar{k} k)]$.

Then $[\Phi_n(ky_{t_{n-1}})](\Phi_n \bar{k}) = [(\Phi_{n-1} k) x_{t_{n-1}}](\Phi_{n-1} \bar{k}) =$

$= [(\Phi_{n-1} y_1) x_{t_{n-1}}]\{[(\Phi_{n-1} k)(\Phi_{n-1} y_1)][(\Phi_{n-1} \bar{k})(\Phi_{n-1} k)]\} =$

$= [(\Phi_{n-1} y_1) x_{t_{n-1}}] \Phi_{n-1}[(ky_1)(\bar{k} k)] = \Phi_n\left[(ky_{t_{n-1}})\bar{k}\right]$. Similarly,

$\Phi_n\{[(y_1 y_{t_{n-1}}) k](\bar{k} y_{t_{n-1}})\} = \Phi_n\left[(y_1 k)(\bar{k} y_1)\right] = \Phi_{n-1}\left[(y_1 k)(\bar{k} y_1)\right] =$

$= [(\Phi_{n-1} y_1)(\Phi_{n-1} k)][(\Phi_{n-1} \bar{k})(\Phi_{n-1} y_1)] =$

$= \{[(\Phi_{n-1} y_1) x_{t_{n-1}}](\Phi_{n-1} k)\}[(\Phi_{n-1} \bar{k}) x_{t_{n-1}}] =$

$= \{\Phi_n\left[(y_1 y_{t_{n-1}}) k\right]\}\left[\Phi_n(\bar{k} y_{t_{n-1}})\right]$. The proofs for the other 14 forms of products on $K_n$ are similar. So every $\Phi_n : K_n \to G_n$ is an isomorphism.

We now define $\Phi : K \to H$ as follows: for $\alpha \in K_n - K_{n-1}$, $\Phi \alpha \equiv \Phi_n \alpha$. Then for any $\{\alpha, \beta\} \subseteq K$, with $\alpha \in K_n - K_{n-1}$ and $\beta \in K_m - K_{m-1}$,

$\Phi(\alpha\beta) = \Phi_M(\alpha\beta) = (\Phi_M \alpha)(\Phi_M \beta) = (\Phi_n \alpha)(\Phi_m \beta) = (\Phi \alpha)(\Phi \beta)$, where $M = \max\{n, m\}$. Using the definition of the $\Phi_n$'s it is straightforward to prove that $\Phi$ is 1-1 and onto $H$. So $H$ and $K$ are isomorphic. ∎

**Corollary 9.** *A countable anti-rectangular $AG$-band is a union of a countable number of disjoint, isomorphic copies of G.*

**Corollary 10.** *A countable anti-rectangular $AG$-band is isomorphic to a proper subgroupoid of itself.*

*Proof*: Consider H in Construction 1. Let $J_1 = \{a, ax_1, x_1a, x_1\}$. For $1 \prec n$ define $J_n$ by induction as $J_n = J_{n-1}^{x_n}$. Then $J = \bigcup_{n=1}^{\infty} J_n$, with the multiplication inherited from *H*, is a proper, countable, anti-rectangular *AG*-sub-band of *H*. By Theorem 8, *J* and *H* are isomorphic. ∎

It follows from Result 9, Corollary 5 and Theorem 8 that:

**Corollary 12.** *If* $R \in ARAGB$, *R is finite or countable and* $R \cong \overline{R}$ *then* $R \cong \overline{R}$.

## 4. Non-anti-rectangular *AG*-groupoids that are an anti-rectangular *AG*-band of anti-rectangular *AG*-bands.

Looking closely at Result 1 it is natural to wonder whether $Y \in ARAGB$ implies $B \in ARAGB$. The converse statement is trivial, since any band $B \in ARAGB$ is an anti-rectangular band $Y = B \in ARAGB$ of trivial anti-rectangular *AG*-bands. However, there is an *AG*-groupoid $\overline{G} \notin ARAGB$ that is an anti-rectangular *AG*-band *Y* of anti-rectangular *AG*-bands $Y_\alpha (\alpha \in Y)$. In fact $\overline{G}$ has order 16, which is the minimal order for such an *AG* band that is *not* anti-rectangular, as we proceed to prove. We also show that $\overline{G}$ is unique up to isomorphism and that any *AG*-band $G^* \notin ARAGB$ that is an anti-rectangular *AG*-band *Y* of anti-rectangular *AG*-bands $Y_\alpha (\alpha \in Y)$ must contain a copy of $\overline{G}$. If $a \in Y_\alpha$ then we will denote $Y_\alpha$ as $Y_a$.

**Lemma 12.** *If* $G^* \in AG$ *is an anti-rectangular AG-band Y of anti-rectangular AG-bands* $Y_\alpha (\alpha \in Y)$ *then*

12.1) $G^*$ is cancellative,
12.2) for any $\{a,b\} \subseteq G^*$, $|Y_a| = |Y_b|$ and
12.3) for any $\{a,b\} \subseteq G^*$, $aba = b$ if and only if $bab = a$.

*Proof*: **12.1)** Suppose that $a \in Y_\alpha = Y_a$, $b \in Y_\beta = Y_b$ and $x \in Y_\gamma = Y_x$. If $xa = xb$ then $\gamma\alpha = \gamma\beta$ and, since *Y* is cancellative, $\alpha = \beta$. Then $aba = b$ and so $bx = (aba)x = (xa)(ab) = (xb)(ab) = [(ab)b]x = (ba)x$. Hence, $(xax)b = (bx)(xa) = [(ba)x](xa) = (xax)(ba)$. But $(b, ba, xax) \subseteq Y_\beta$ and $Y_\beta$ is cancellative. Therefore $b = ba = bb$ and so $a = b$. Dually, if $ax = bx$ then $a = b$. Therefore $G^*$ is cancellative.

**12.2)** Now let $x \in Y_\alpha = Y_a$. Then $(ab)x \in Y_\beta$. Since $Y_\beta$ is cancellative $|Y_\alpha| \leq |Y_\beta|$. Dually $|Y_\beta| \leq |Y_\alpha|$ and so $|Y_\alpha| = |Y_\beta|$.

**12.3)** If $aba = b$ then $bab = a[(bab)a] = a[(ba)(aba)] = a[(ba)b]$. But $\{a, bab\} \subset Y_\alpha$ and $Y_\alpha$ is cancellative. Hence $a = bab$. Dually, $bab = a$ implies $aba = b$. ∎

Now suppose that $G^* \in \mathbf{AG}$ is an anti-rectangular $\mathbf{AG}$-band $Y$ of anti-rectangular $\mathbf{AG}$-bands $Y_\alpha (\alpha \in Y)$. If $G^*$ is <u>not</u> anti-rectangular then by Lemma 12 we have $\{a, b, x, y\} \subseteq G^*$ with $aba = y \neq b$, $bab = x \neq a$, $\{a, x, ax, xa\} \subseteq Y_a$ and $\{b, y, by, yb\} \subseteq Y_b$.

Since $G^*$ is an anti-rectangular $\mathbf{AG}$-band $Y$ of anti-rectangular $\mathbf{AG}$-bands $Y_\alpha (\alpha \in Y)$ and from remarks in the paragraph preceding Example 1 it follows that $\{a,, x, ax, xa\} \equiv G_a$, $\{b, y, by, yb\} \equiv G_b$, $\{ab, xy, (ab)(xy), (xy)(ab)\} \equiv G_{ab}$ and $\{ba, yx, (ba)(yx), (yx)(ba)\} \equiv G_{ba}$ are disjoint, isomorphic copies of $G$ contained in $Y_a, Y_b, Y_{ab}$ and $Y_{ba}$ respectively. We proceed to demonstrate that the union $\overline{G} = \bigcup G_g \, (g \in G)$ of these four copies of $G$ is a subgroupoid of $G^*$ and an anti-rectangular $\mathbf{AG}$-band $G$ of anti-rectangular $\mathbf{AG}$-bands $G_g \, (g \in \{a, b, ab, ba\})$.

Recall that $G^*$ is a cancellative $\mathbf{AG}$-band. We have $y = (ab)a$. Then
$(ab)x = (xb)a = [(bab)b]a = [b(ba)]a = (aba)b = yb$,
$(ab)(ax) = (aba)[(ab)x] = y(yb) = by$ and
$(ab)(xa) = [(ab)x](aba) = (yb)y = b$.

We have shown that $G_b = \{b, y, by, yb\} = (ab)\{xa, a, ax, x\} = (ab)G_a$. Similarly we can calculate that
$G_{ab} = \{ab, xy, (ab)(xy), (xy)(ab)\} = \{a, ax, xa, x\}b = G_a b$ and

$G_{ba} = \{ba, yx, (ba)(yx), (yx)(ba)\} = b\{a, xa, x, ax\} = bG_a$. We can then calculate the Cayley table consisting of the 256 products of pairs of elements of $\overline{G}$. In order to have sufficient space to show the Cayley table we define the following two ordered 16-tuples as equal: $(1,2,3,4,5,6,7,8,9,10,11,12,13,14,15,16) =$
$(a, x, ax, xa, b, y, by, yb, ab, xy, (ab)(xy), (xy)(ab), ba, yx, (ba)(yx), (yx)(ba))$.

The Cayley table of $\overline{G}$ is shown below. Table 4 is derived from Table 3.

**Table 3**

| $\overline{G}$ | 1 | 2 | 3 | 4 | 5 | 6 | 7 | 8 | 9 | 10 | 11 | 12 | 13 | 14 | 15 | 16 |
|---|---|---|---|---|---|---|---|---|---|---|---|---|---|---|---|---|
| 1 | 1 | 3 | 4 | 2 | 9 | 11 | 12 | 10 | 16 | 14 | 13 | 15 | 6 | 8 | 7 | 5 |
| 2 | 4 | 2 | 1 | 3 | 12 | 10 | 9 | 11 | 13 | 15 | 16 | 14 | 7 | 5 | 6 | 8 |
| 3 | 2 | 4 | 3 | 1 | 10 | 12 | 11 | 9 | 15 | 13 | 14 | 16 | 5 | 7 | 8 | 6 |
| 4 | 3 | 1 | 2 | 4 | 11 | 9 | 10 | 12 | 14 | 16 | 12 | 13 | 8 | 6 | 5 | 7 |
| 5 | 13 | 15 | 16 | 14 | 5 | 7 | 8 | 6 | 2 | 4 | 3 | 1 | 12 | 10 | 9 | 11 |
| 6 | 16 | 14 | 13 | 15 | 8 | 6 | 5 | 7 | 3 | 1 | 2 | 4 | 9 | 11 | 12 | 10 |
| 7 | 14 | 16 | 15 | 13 | 6 | 8 | 7 | 5 | 1 | 3 | 4 | 2 | 11 | 9 | 10 | 12 |
| 8 | 15 | 13 | 14 | 16 | 7 | 5 | 6 | 8 | 4 | 2 | 1 | 3 | 10 | 12 | 11 | 9 |
| 9 | 6 | 8 | 7 | 5 | 13 | 15 | 16 | 14 | 9 | 11 | 12 | 10 | 4 | 2 | 1 | 3 |
| 10 | 7 | 5 | 6 | 8 | 16 | 14 | 13 | 15 | 12 | 10 | 9 | 11 | 1 | 3 | 4 | 2 |
| 11 | 5 | 7 | 8 | 6 | 14 | 16 | 15 | 13 | 10 | 12 | 11 | 9 | 3 | 1 | 2 | 4 |
| 12 | 8 | 6 | 5 | 7 | 15 | 13 | 14 | 16 | 11 | 9 | 10 | 12 | 2 | 4 | 3 | 1 |
| 13 | 9 | 11 | 12 | 10 | 2 | 4 | 3 | 1 | 8 | 6 | 5 | 7 | 13 | 15 | 16 | 14 |
| 14 | 12 | 10 | 9 | 11 | 3 | 1 | 2 | 4 | 5 | 7 | 8 | 6 | 16 | 14 | 13 | 15 |
| 15 | 10 | 12 | 11 | 9 | 1 | 3 | 4 | 2 | 7 | 5 | 6 | 8 | 14 | 16 | 15 | 13 |
| 16 | 11 | 9 | 10 | 12 | 4 | 2 | 1 | 3 | 6 | 8 | 7 | 5 | 15 | 13 | 14 | 16 |

**Table 4.** (for $g, \overline{g} \in G_a = \{a, x, ax, xa\}$)

| $\overline{G}$ | $\overline{g}$ | $(ab)\overline{g}$ | $\overline{g}b$ | $b\overline{g}$ |
|---|---|---|---|---|
| $g$ | $g\overline{g}$ | $\{x[g(a\overline{g})]\}b$ | $b\{[a(\overline{g}g)]x\}$ | $(ab)[(ga)\overline{g}]$ |
| $(ab)g$ | $b[(xa)(\overline{g}g)]$ | $(ab)(g\overline{g})$ | $(xg)(\overline{g}a)$ | $[(g\overline{g})a]b$ |
| $gb$ | $(ab)[(\overline{g}a)(g\overline{g})]$ | $b[(\overline{g}g)(xa)]$ | $(g\overline{g})b$ | $\overline{g}[(ag)x]$ |
| $bg$ | $(\overline{g}g)b$ | $\overline{g}(gx)$ | $(ab)[g(x\overline{g})]$ | $b(g\overline{g})$ |

Notice that Table 4 yields the following Cayley table in set theoretic notation:

**Table 5.**

| $\overline{G}$ | $G_a$ | $G_b = (ab)G_a$ | $G_{ab} = G_a b$ | $G_{ba} = bG_a$ |
|---|---|---|---|---|
| $G_a$ | $G_a$ | $G_{ab}$ | $G_{ba}$ | $G_b$ |
| $G_b = (ab)G_a$ | $G_{ba}$ | $G_b$ | $G_a$ | $G_{ab}$ |
| $G_{ab} = G_a b$ | $G_b$ | $G_{ba}$ | $G_{ab}$ | $G_a$ |
| $G_{ba} = bG_a$ | $G_{ab}$ | $G_a$ | $G_b$ | $G_{ba}$ |

Note that the subscripts of the $G_g$'s $(g \in G)$ multiply in exactly the same way as the elements of $G$ in Example 1. The fact that $\overline{G} \in AG$ follows from the fact that $\overline{G}$ is a subgroupoid of $G^* \in ARAGB \subseteq AG$. This proves that $\overline{G}$ is an anti-rectangular $AG$-band $G$ of anti-rectangular $AG$-bands $G_g$ $(g \in \{a, b, ab, ba\})$. Note however that $\{a, b, ab, ba\}$ is not even a subgroupoid of $\overline{G}$!

We have therefore proved the following:

**Theorem 12.** $\overline{G} \in AG$ and is an anti-rectangular $AG$-band $G$ of (four) isomorphic copies of $G$. However $\overline{G} \notin ARAGB$. Also, if $G^* \in AG$, $G^*$ is an anti-rectangular $AG$-band $Y$ of anti-rectangular $AG$-bands $Y_\alpha (\alpha \in Y)$ and $G^* \notin ARAGB$ then $G^*$ contains an isomorphic copy of $\overline{G}$. Finally, $\overline{G}$ is the smallest $AG$-groupoid $K$ satisfying the conditions $K \notin ARAGB$ and $K$ is an anti-rectangular band $Y$ of anti-rectangular $AG$-bands $Y_\alpha (\alpha \in Y)$.